\documentclass[12pt]{article}
\usepackage{amsmath, amsthm, amssymb, graphicx}

\newtheorem{theorem}{Theorem}[section]
\newtheorem{lemma}[theorem]{Lemma}
\theoremstyle{definition}
\newtheorem{definition}[theorem]{Definition}

\theoremstyle{remark}

\numberwithin{equation}{section}

\begin{document}
\title{Monotone Periodic Orbits for Torus Homeomorphisms}
\author{Kamlesh Parwani}
\date{September 1, 2003.}
\maketitle
\begin{abstract}
Let $f$ be a homeomorphism of the torus isotopic to the identity and suppose
that there exists a periodic orbit with a non-zero rotation vector $(\frac{p}{q},\frac{r}{q})$, then $f$ has a topologically monotone periodic orbit with the same rotation vector.
\end{abstract}

\section*{Introduction}

In this article we prove a theorem about the existence of topologically
monotone periodic orbits on the torus. The concept of monotone orbits on the
annulus is certainly not new; it goes back to Aubry and Mather's proof of
the existence of orbits whose radial order is preserved by an
area-preserving twist map of the annulus. These orbits, that have their
radial order preserved by the map, are called Birkhoff orbits (see \cite
{Katok}) or monotone orbits.

This notion of monotone orbits inspired the definition of topologically monotone
orbits in \cite{Boyland}, where Boyland proved that any homeomorphism of the
annulus isotopic to the identity that has a periodic orbit with a non-zero rotation
number $\frac{p}{q}$ also has a topologically monotone periodic orbit with
the same rotation number. A topologically monotone periodic orbit has the
property that the isotopy class of the map, keeping the periodic orbit fixed
as a set, is of finite order. The main tool used in Boyland's proof is
Nielsen-Thurston theory.

In \cite{Llibre&Mackay}, Llibre and Mackay asked whether a similar result
was true for torus homeomorphisms. The goal of this paper is to answer that
question by proving the same theorem on the torus.

\textbf{Main Theorem.} \textit{If }$f$\textit{\ is a torus homeomorphism
isotopic to the identity that has a periodic orbit with a non-zero rotation vector }$%
\left( \frac{p}{q},\frac{r}{q}\right) $\textit{, then }$f$\textit{\ also has
a topologically monotone periodic orbit with the same rotation vector.}

\smallskip

There are some immediate complications one encounters while trying to
generalize the theorem to the torus. First, the torus has rotation vectors
instead of rotation numbers. Then on the annulus, under certain restrictions
to the rotation number, we can only get the $pA$ (pseudoAnosov) isotopy
class or the finite order isotopy class, but for the torus, there is the
reducible isotopy class to deal with also. These concepts will be introduced
in Section 1, and then in Section 2 we prove the main theorem of this paper.

It should be noted that LeCalvez has proved the existence of topologically
monotone periodic orbits on the torus under the assumption that the maps are
smooth by using variational techniques (see \cite{LeCalvez}). Since we are
dealing with homeomorphisms, we rely solely on topological methods.

\section{Definitions and important results}

\subsection{Rotation vectors}

Let $f$ be a homeomorphism of the torus which is isotopic to the identity
and let $F$ be its lift to the universal cover $\widetilde{T^{2}}$, the
plane. Let $()_{1}$ and $()_{2}$ be the projections of a point in the plane
to the $x$-axis and the $y$-axis respectively and let $x$ be a point on the
torus $T^{2}$ with $\widetilde{x}$ as its lift. Then the \textbf{rotation
vector} of $\widetilde{x}$, with respect to a lift $F$, is defined as
following if the limit exists.

\begin{equation*}
\rho (\widetilde{x},F)=\left( \lim_{n\to \infty }\left( \frac{F^{n}(%
\widetilde{x})-\widetilde{x}}{n}\right) _{1},\lim_{n\to \infty }\left( \frac{%
F^{n}(\widetilde{x})-\widetilde{x}}{n}\right) _{2}\right)
\end{equation*}

If $x$ is a periodic point, say of period $q$, then the rotation vector is
always well defined and can be written as $(\frac{p}{q},\frac{r}{q})$ for
some integers $p$ and $r$. In fact, the rotation vector for any point on the
orbit of $x$ is the same, and so, we can associate the vector $(\frac{p}{q},%
\frac{r}{q})$ to the periodic orbit. Periodic orbits with rotation vector
$(\frac{p}{q},\frac{r}{q})$ and least period $q$ will be called $(p,r,q)$ orbits. Note that a
$(p,r,q)$ orbit has the same rotation vector as a $(pt,rt,qt)$ orbit, where $t$
is some positive integer.

The covering space of the torus comes naturally equipped with two important
covering translations. Define $X(\widetilde{x})=\widetilde{x}+(1,0)$ and $Y(%
\widetilde{x})=\widetilde{x}+(0,1)$. Clearly, the rotation vector of a point
depends on the lift, and the relationship is $\rho (\widetilde{x}%
,Y^{m}X^{n}F)=\rho (\widetilde{x},F)+(n,m)$, where $n$ and $m$ are integers.
So when we discuss periodic orbits with a certain rotation vector in this
paper, we assume the existence of some lift for which that rotation vector
is realized.  Also, when we start with a periodic orbit, say $(p,r,q)$ orbit,
and then prove that another $(p,r,q)$ orbit exists, it is to be understood
that both rotation vectors are calculated with respect to the same lift.

\subsection{The Nielsen-Thurston classification theorem and braids}

Every orientation preserving homeomorphism of an orientable surface with
negative Euler characteristic is isotopic to a homeomorphism $g$ such that
either

a) $g$ is finite order, or

b) $g$ is pseudoAnosov ($pA$), or

c) $g$ is reducible.

A map $g$ is said to be \textit{reducible} if there is a disjoint collection 
$C$ of non-parallel, non-peripheral simple disjoint curves such that $g$
leaves invariant the union of disjoint regular neighborhoods of curves in $C$%
, and the first return map on each complementary component is either of
finite order or $pA$.

This classification theorem was first announced in \cite{Thurston} and the
proofs appeared later in \cite{FLP} and \cite{Casson}.

\smallskip

The torus doesn't have negative Euler characteristic but, following Handel
as in \cite{Handel}, we will examine the isotopy class relative to a
periodic orbit; this will introduce punctures and provide the negative Euler
characteristic to apply the Nielsen-Thurston Classification Theorem. When
the isotopy class relative to a given periodic orbit is of finite order, the
periodic orbit is called a \textit{finite order periodic orbit}, and \textit{%
reducible} and $pA$ \textit{periodic orbits} are defined similarly. The
isotopy class relative to a periodic orbit is also referred to as the 
\textit{braid} of the periodic orbit.

\begin{definition}
Let $x$\ and $y$\ be two distinct periodic points of least period $n$\ for
homeomorphisms $f$ and $g$ respectively of the same orientable surface $S$.
Then the orbit of $x$\ ($O(x)$) and the orbit of $y$\ ($O(y)$) have the same 
\textit{braid} if there exists an orientation-preserving homeomorphism $h$\
of $S$ with the property that $h$\ maps $O(x)$\ onto $O(y)$ and the isotopy
class of $h^{-1}fh$\ relative to the orbit of $y$\ is the same as the
isotopy class of $g$\ relative to the orbit of $y$, that is, $%
[h^{-1}fh]_{O(y)}=[f]_{O(y)}$.
\end{definition}

A periodic orbit has a \textbf{trivial} braid if the isotopy class relative to the
periodic orbit is of finite order, that is,
there exists a homeomorphism $g$ isotopic to $f$, relative to the periodic
orbit, such that $g^{n}=identity$ for some $n$.  In other words, finite order periodic orbits
have trivial braids.  These periodic orbits are considered to be topologically monotone.

A periodic orbit has a \textbf{non-trivial} braid 
if the isotopy class relative to the periodic orbit is not of finite order. In other words,
periodic orbits with non-trivial braids are either reducible periodic orbits or are $pA$ periodic orbits.
These periodic orbits are not topologically monotone.

Boyland defined a natural partial order $(\vartriangleright)$ into these
braids. If $\alpha $ and $\beta $ are two braids of periodic orbits, then $%
\alpha \vartriangleright \beta $ if and only if the existence of a periodic orbit with braid $\alpha $
in any homeomorphism $f$ on a given surface implies the existence of a periodic orbit with braid $\beta $ for the
same $f$. The proof of the fact that this is an actual partial order is not
easy and is in \cite{Boyland}. He also proved that a $pA$ periodic orbit is
strictly above (in the partial ordering) all other periodic orbits that are
present in the $pA$ representative of the isotopy class relative to the $pA$ periodic orbit.

The existence of topologically monotone periodic orbits on the annulus in 
\cite{Boyland} is established by showing that periodic orbits with non-trivial braids
force the existence of periodic orbits with trivial braids and the
same rotation number (non-trivial $\vartriangleright $ trivial).

\begin{theorem}[Boyland]
Let f be a homeomorphism of the annulus isotopic to the identity.  If f has a
periodic orbit with non-zero rotation number $\frac{p}{q}$, then f also has a topologically monotone
periodic orbit with the same rotation number.
\end{theorem}

Essentially, we follow the same strategy on the torus and the proof of the
main theorem in Section 2 is broken into two parts, reducible $%
\vartriangleright $ finite order and $pA$ $\vartriangleright $ finite order. 

We will also need the following result which can be obtained from the arguments
in \cite{Hall&Boyl} and is also proved in \cite{thesis}.

\begin{theorem}
Let f be a homeomorphism of the torus that is isotopic to the identity and
has a pA periodic orbit. Let g be the pA representative of the isotopy class relative to this
orbit and let G be its lift to the plane that fixes the lifts of all the
points in the pA orbit. Then G has a dense orbit.
\end{theorem}

We will use this theorem in the next section to prove $pA$ $\vartriangleright 
$ trivial. 

\section{Finite order periodic orbits on the torus}

In this section we prove the main theorem by showing that finite order
periodic orbits are on the bottom in the partial ordering of periodic orbits
for torus homeomorphism isotopic to the identity, that is, the reducible $%
\vartriangleright $ finite order and $pA\vartriangleright $ finite order. We
restrict our attention to periodic orbits with the same non-zero rotation vector, say 
$(\frac{p}{q},\frac{r}{q})$, and assume that there are no common factors between $p$, $r$, and 
$q$. Later we reduce the general case, in which there may be a common factor
between $p$, $r$, and $q$, to the case of no factors.
\begin{theorem}
Let $f:T^{2}\rightarrow T^{2}$ be a homeomorphism isotopic to
the identity. Suppose there exists a $(p,r,q)$ orbit such that $gcd(p,r,q)=1$,
then there exists a finite order $(p,r,q)$ orbit.
\end{theorem}
\begin{proof}
If the given periodic orbit is already of finite order type,
then we're done. If not, the proof breaks down in to two cases---the
periodic orbit is either of reducible type or it is of $pA$ type. These are
handled separately below.
\smallskip

\textbf{Reducible Case}.
\smallskip

We assume that we obtain a reducible isotopy class, keeping the periodic
orbit fixed. A reducing curve can be of two types---essential or
non-essential (this is with respect to the unpunctured torus). A reducing
curve cannot be non-essential for this would imply that there is a common
factor between $p,r,q$. In fact, the common factor would be exactly the
number of punctures contained in the disc bounded by the non-essential curve.

So any reducing curve must be essential. All the essential reducing curves are disjoint, 
and thus, split the torus into parallel
annuli. These annuli have an equal number of punctures, say $n$ punctures,
in their interiors and are permuted by the action of the map $g$, which
represents this reducible isotopy class. 

In this case, we will show that we can obtain a finite order isotopy class,
keeping a periodic orbit with the same period and rotation vector (it may
not be the same periodic orbit that we start with). Let $%
A_{0},A_{1},A_{2},...,A_{m-1}$ be the annuli that $g$ permutes, numbered so
that $g(A_{k})=A_{k+1}$ mod $m$ and $nm=q$.

\begin{figure}[ht]
\centerline{\includegraphics[height=2.0in]{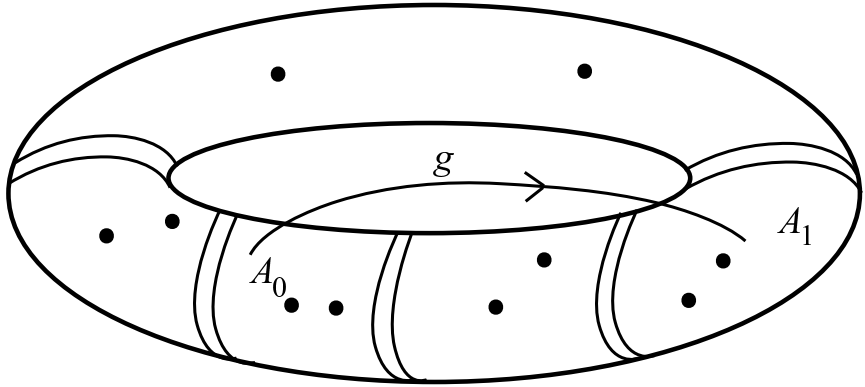}}
\caption{The Reducible Case}
\end{figure}

CASE 1.

Suppose that the maps $g^{m}:A_{k}\rightarrow A_{k}$ relative to the
punctures are all of finite order. It is easy to see that all these maps are
conjugate to each other, and so if one is of finite order, then all are of
finite order. Since all finite order maps on the annulus are conjugate to rotations, it
follows that $g^{mn}$ (or $g^{q}$) is the identity in each annulus.

We will now argue that this implies that $g^{q}$ is isotopic to the identity on
the entire torus relative to the $(p,r,q)$ orbit.  The complements of the interiors of annuli
containing the punctures (the $A_{k}$'s) are closed annuli that do not contain any punctures.
Observe that $g^{q}$ fixes the boundary components of these unpunctured annuli since it fixes the
boundary components of the $A_{k}$'s.  So if $g^{q}$ is isotopic to the identity in all these unpunctured annuli relative
to their boundaries (keeping the boundaries fixed throughout the isotopy) then $g^{q}$
is isotopic to the identity on the entire torus relative to the $(p,r,q)$ orbit, because we already know
that $g^{q}$ is the identity on the annuli containing the punctures.  Now suppose $g^{q}$ is not isotopic
to the identity in one of these unpunctured annuli relative to its boundary components, then $g^{q}$ must be isotopic to some non-trivial Dehn twists.  It is easy
to see that the maps on all these annuli are all conjugate to each other so $g^{q}$ is isotopic to the same
Dehn twists in each unpunctured annulus.  However, $g^{q}$ is isotopic to identity on the entire torus
when the punctures are allowed to move, because $g$ is isotopic to $f$ and $f$ is isotopic to the identity
by assumption.  If we have non-trivial Dehn twists that don't cancel each other out (because they are identical),
$g^{q}$ is not isotopic to identity, which is a
contradiction.  It follows that $g^{q}$ is isotopic to the identity on
the entire torus relative to the $(p,r,q)$ orbit.

To show that $(p,r,q)$ orbit is topologically monotone, we require a map isotopic to $g$
relative to the orbit that is of finite order. Such a map is guaranteed
by Fenchel's solution to the Nielsen Realization problem for finite solvable groups (see Chapter 3 in
\cite{Zieschang}) which provides a map $h$ isotopic to $g$, relative to the $
(p,r,q)$ orbit, such that $h^{q}$ is the identity. This shows that the $
(p,r,q)$ orbit is of finite order type, that is, it is topologically
monotone.

CASE 2.

Suppose that the maps $g^{m}:A_{k}\rightarrow A_{k}$ relative to the
punctures are all $pA$. It now follows from Boyland's proof of Theorem 1.2 in
\cite{Boyland} that there exists a finite order periodic orbit with the same period and
rotation number in each $A_{k}$.
Since all these annulus maps are conjugate, the periodic orbits connect in
the torus to give a periodic orbit with the same rotation vector and period as the
originally punctured orbit. This reduces to Case 1 and we can find an
isotopy relative to the new orbit such that the isotoped map is of finite
order. Because the periodic orbits in the $pA$ components are unremovable
(see \cite{Boyl:stability}), this periodic orbit existed in the original map 
$f$.

We have actually established a stronger result. If we do obtain a finite order
periodic orbit which is distinct from the one we began with, then it is
strictly below the original orbit in the partial order. This is because the
only way for the original orbit to not be of a finite order type is for the
reducible components to be $pA$, that is, the maps $g^{m}:A_{k}\rightarrow
A_{k}$ are $pA$. And $pA$ orbits are strictly above all periodic orbits that
are present in the $pA$ representative of the isotopy class (see \cite{Boyland}).

\smallskip

\textbf{\textit{pA} Case}.

In Lemma 2.3, we will prove that any $pA$ $(p,r,q)$ orbit forces another $(p,r,q)$ orbit.
Boyland proved that this other periodic orbit is strictly below the $pA$ orbit in the partial order (see 
\cite{Boyland}). Furthermore, there are only finitely many periodic orbits
of any given period in any $pA$ map. So consider a minimal $(p,r,q)$
orbit in the partial order. If it's $pA$, there is another $(p,r,q)$ orbit below it---so
it's not minimal. If it's reducible and not of finite order, then it forces
another $(p,r,q)$ finite order orbit (by the argument above for the reducible case).
Thus, any minimal $(p,r,q)$ orbit must be of finite order and there is at least one minimal orbit.
\end{proof}

It now remains to prove Lemma 2.3. We shall appeal to the following result
in \cite{Fathi} for the existence of a fixed point of positive index. The
proof is based on the ideas used to demonstrate the Brouwer Plane
Translation Theorem (see \cite{Franks}).

\begin{theorem}[Fathi]
Let $G:R^{2}\rightarrow R^{2}$\ be an orientation preserving homeomorphism
which possesses a non-wandering point, then G has a fixed point. If G has
only isolated fixed points, then it has a fixed point of positive index.
\end{theorem}

\begin{lemma}
Let $g:T^{2}\rightarrow T^{2}$\ be the $\mathit{pA}$\
representative obtained from $\mathit{f}$ keeping the $\mathit{(p,r,q)}$
orbit fixed throughout the isotopy. Then $\mathit{g}$\ has another $\mathit{%
(p,r,q)}$\ periodic orbit.
\end{lemma}
\begin{proof}
Let $g$ be the $pA$ representative obtained relative to the $(p,r,q)$
orbit. Let $G$ be the lift
to the plane that realizes the rotation vector $(\frac{p}{q},\frac{r}{q})$ for the $pA$ orbit
and then consider $X^{-p}Y^{-r}G^{q}$; call this map $H$. Then by Theorem 1.3, we know that $H$ has a dense
orbit and so there is a non-wandering point. Since the periodic points (of
any given period) are isolated in a $pA$ map (see \cite {FLP}), the fixed points of $H$ on the
plane are also isolated. So, by Theorem 2.2, we have a fixed point with
positive index for $H$, which is a $(p,r,q)$ orbit for $g$, and
each point in the orbit is fixed with positive index for $g^{q}$.

The $pA$ $(p,r,q)$ orbit is the location of all the one-prongs or needles
and these needles have index zero. Since there is a fixed point for $g^{q}$ with
positive index and because the indices have to add up to zero, we have
actually shown that the $pA$ periodic orbit forces at least two other
$(p,r,q)$ periodic orbits.
\end{proof}

\begin{proof}[Proof of the Main Theorem]
Let $f$ be a homeomorphism of the torus
isotopic to the identity and suppose there exists a periodic orbit with
rotation vector $(\frac{p}{q},\frac{r}{q})$. Also assume that this is a $(pt,rt,qt)$
periodic orbit, where $t$ is a positive integer and there are no common factors between $p$, $r$, and $q$. Let $F$
be the lift to the plane and consider $X^{-p}Y^{-r}F^{q}$. We obtain a
periodic point of period $t$ for $X^{-p}Y^{-r}F^{q}$. Then by Theorem 2.2,
we also obtain a fixed point for $X^{-p}Y^{-r}F^{q}$, which has period $q$
for $f$ and rotation vector $(\frac{p}{q},\frac{r}{q})$, that is, it is a $(p,r,q)$ orbit.
Thus, without loss of generality, we may assume that we have a $(p,r,q)$ orbit where there are no common
factors between $p$, $r$, and $q$. Now, by Theorem 2.1, we obtain a $(p,r,q)$ orbit that is
topologically monotone and it has the desired rotation vector $(\frac{p}{q},\frac{r}{q})$.
\end{proof}

It is natural to ask if there is a similar theorem for non-periodic orbits
with irrational rotation vectors. This question is unanswered even for the
annulus. A similar theorem about the existence of monotone periodic orbits
has been proved for periodic orbits on surfaces of higher genus (see \cite
{Parwani}).

\section*{Acknowledgments}

The author would like to thank John Franks and Philip Boyland for several
useful and stimulating conversations.

\end{document}